# Approximate solutions by truncated Taylor series expansions of nonlinear differential equations and the related shadowing property

M. De la Sen

Instituto of Research and Development of Processes. University of the Basque Country

Campus of Leioa (Bizkaia) - P.O. Box 644- Bilbao; Barrio Sarriena, 48940- Leioa. SPAIN

**Abstract**. This paper investigates the errors of the solutions as well as the shadowing property of a class of nonlinear differential equations which possess unique solutions on a certain interval for any admissible initial conditions. The class of differential equations are assumed to be approximated by well-posed truncated Taylor series expansions up to a certain order obtained about certain, in general non-periodic, sampling points $t_i \in [t_0, t_J]$ for $i = 0,1,..., J$ of the solution.

## 1. Introduction

This paper investigates the errors of the solutions of nonlinear differential equations $\dot{y}(t) = f(y(t), t)$, where $f \in C^{(n+1)}(\mathbf{R}^n \times (t_0, t_J); \mathbf{R}^n)$, provided they exist and are unique for each given admissible initial condition $y(t_0) = y_0$, with respect to the solutions of its approximate differential equations $\dot{x}(t) = \sum_{k=0}^{\ell} \frac{f^{(k)}(x(t_i), t_i)}{k!}(x(t) - x(t_i))^k$ ; $x(t_0) = x_0$, for any given nonnegative integer $\ell \leq n$, obtained from truncated Taylor expansions of the solutions about certain sampling points $t_i \in [t_0, t_J]$ for $i = 0,1,..., J$. It is assumed that if a unique solution exists on some interval $[t_0, t_J]$ and that the choice of the sampling points is such that the inter-sample intervals [7-10] are subject to a certain maximum allowable upper-bound then the error of the solution in the whole interval $[t_0, t_J]$ satisfies a prescribed norm bound. Using the obtained results, the shadowing property [1-6] of the true solution with respect to the approximate one is investigated in the sense that "shadowing" initial conditions of the true solution exist, for each initial condition of the approximate differential equation, such that any approximated solution trajectory on the interval of interest is arbitrarily close to the true one under prescribed allowable maximum norms of the error between both the true and the approximate solutions. The problem is extended to the case when the approximated solution is perturbed either by a sequence of a certain allowable size at the sampling points or with perturbation functions of a certain size in norm about the whole considered interval. The main tool involved in the analysis is an "ad hoc" use of a known preparatory theorem to the celebrated Bernstein´s theorem, [11], which gives an upper-bound for the maximum norm of the error in-between both the true and the approximate solutions. The results are extendable functional equations involving the presence of delays. [12-13].

## 2. Calculation of the exact solution from Taylor series expansion

**Lemma 2.1**. Assume that $f \in C^{(n+1)}((a,b); \mathbf{R}^n)$ and divide the real interval $(a,b)$ into $J$ subintervals with points $y_n \in [a,b]$ such that

$$a \equiv y_0 < y_1 < y_2 < .... < y_{J-1} < y_J \equiv b.$$ Then:

$$y_{i+1} = y_0 + \int_0^{\bar{h}_i} f(y)dy = y_0 + \sum_{j=-1}^{i} \int_{\bar{h}_j}^{\bar{h}_{j+1}} f(y)dy = y_0 + \sum_{j=-1}^{i} \int_0^{h_j} f(y + \bar{h}_j)dy$$

$$= y_0 + \sum_{j=-1}^{i} \sum_{k=0}^{n} \frac{f^{(k)}(c_j)}{k!} \int_0^{h_j} (y + \bar{h}_j - c_{j+1})^k dy + \frac{1}{n!} \sum_{j=-1}^{i} \int_0^{h_j} \int_{c_{j+1}}^{y_{j+1}} (y + \bar{h}_j - t)^n f^{(n+1)}(t) dt dy$$

(2.1)

where $h_n = y_{n+1} - y_n$, $\bar{h}_n = y_{n+1} - y_0 = \sum_{i=0}^{n} h_i$ for $n = 0,1,....,J-1$ with $\bar{h}_{-1} = 0$, and

$$y_{i+1}(\tilde{h}_{i+1}) = y_0 + \sum_{j=-1}^{i} \sum_{k=0}^{n} \frac{f^{(k)}(c_j)}{k!} \int_0^{h_j} (y + \bar{h}_j - c_{j+1})^k dy + \frac{1}{n!} \sum_{j=-1}^{i} \int_0^{h_j} \int_{c_{j+1}}^{y_{j+1}} (y + \bar{h}_j - t)^n f^{(n+1)}(t) dt dy$$

$$+ \sum_{k=0}^{n} \frac{f^{(k)}(c_{i+1})}{k!} \int_0^{\tilde{h}_{i+1}} (y + \bar{h}_i + \tilde{h}_{i+1} - c_{i+2})^k dy + \frac{1}{n!} \int_0^{\tilde{h}_{i+1}} \int_{c_{j+1}}^{y_{j+1}} (y + \bar{h}_i + \tilde{h}_{i+1} - t)^n f^{(n+1)}(t) dt dy \quad (2.2)$$

; $\forall y \in [y_i, y_{i+1}]$, and any real $c_i \in [y_i, y_{i+1}]$ for $i = 0,1,....,J-1$; $\forall \tilde{h}_i \in [0, h_i]$ for $i = 0,1,....,J-2$.

*Proof*: It follows from a well-known preparatory theorem to Bernstein's theorem [1] since

$$f(y) = \sum_{k=0}^{n} \frac{f^{(k)}(c_i)}{k!} (y - c_i)^k + \frac{1}{n!} \int_{c_i}^{y} (y-t)^n f^{(n+1)}(t) dt \tag{2.3}$$

□

Now, consider the nonlinear ordinary differential equation:

$$\dot{y}(t) = f(y(t), t) \; ; \; y(t_0) = y_0 \tag{2.4}$$

in the real interval $\mathbf{R}^n \times [t_0, t_J]$ such that $f \in C^{(n+1)}(\mathbf{R}^n \times (t_0, t_J); \mathbf{R}^n)$ is Lipschitz-continuous in $[y(t_0) - \theta_0, y(t_0) + \theta_0] \times [t_0, t_J]$. The following result follows from Lemma 2.1:

**Theorem 2.2**. The unique solution of (2.4) in $[t_0, t_J]$ is given by:

$$y(t) = y(t_0) + \sum_{j=0}^{i-2} \int_0^{h_j} \left[ \sum_{k=0}^{\ell} \frac{f^{(k)}(y(t_j), t_j)}{k!} (y(\tau + t_j) - y(t_j))^k + \frac{1}{\ell!} \int_0^{\tau} (y(\sigma + t_j) - y(t_j))^{\ell} f^{(\ell+1)}(y(\sigma + t_j), \sigma + t_j) d\sigma \right] d\tau$$

$$+ \int_{t_{i-1}}^{t} \left[ \sum_{k=0}^{\ell} \frac{f^{(k)}(y(t_j), t_j)}{k!} (y(\tau + t_{i-1}) - y(t_{i-1}))^k + \frac{1}{\ell!} \int_0^{\tau} (y(\sigma + t_{i-1}) - y(t_{i-1}))^{\ell} f^{(\ell+1)}(y(\sigma + t_{i-1}), \sigma + t_{i-1}) d\sigma \right] d\tau$$

(2.5)

; $\forall t \in [t_{i-1}, t_i]$ ; $\forall i \in \bar{J} = \{1, 2, ..., J\}$, $\forall \ell (\in \mathbf{Z}_{0+}) \leq n$, where $t_i \in [t_0, t_J]$ are any arbitrary strictly ordered points such that $t_0 < t_1 < t_2 < .... < t_{J-1} < t_J$ with $h_i = t_{i+1} - t_i$ for $i = 0, 1, ...., J-1$.

*Proof*: Since $f \in C^{(n+1)}(\mathbf{R}^n \times (t_0, t_J); \mathbf{R}^n)$ is Lipschitz-continuous in $[y(t_0) - \theta_0, y(t_0) + \theta_0] \times [t_0, t_J]$ so that the solution $y(t)$ on $[t_0, t_J]$ is unique, provided that $t_J = t_J(\theta_0, t_0)$ for the given $t_0 \in \mathbf{R}$ and some $\theta_0 \in \mathbf{R}^n$ is such that $y(t) \in [y(t_0) - \theta_0, y(t_0) + \theta_0]$ ; $\forall t \in [t_0, t_J]$ and $t_J \in (t_0, t_{\bar{J}}]$ since $f : [y(t_0) - \theta_0, y(t_0) + \theta_0] \times [t_0, t_J] \to \mathbf{R}^n$ is local Lipschitz-continuous as a result. Such a unique solution is given by:

$$y(t) = y_a + \int_a^t f(y(\tau), \tau) d\tau \; ; \; \forall t \in [t_0, t_J] \tag{2.6}$$

Take any set of $J$ strictly ordered points $t_n \in [t_0, t_J]$ satisfying $t_0 < t_1 < t_2 < .... < t_{J-1} < t_J$ with $h_i = t_{i+1} - t_i$ for $i = 0, 1, ...., J-1$, so that:

$$y(t) = y(t_{i-1}) + \int_{t_{i-1}}^t f(y(\tau), \tau) d\tau = y(t_0) + \sum_{j=0}^{i-2} \int_{t_j}^{t_{j+1}} f(y(\tau), \tau) d\tau + \int_{t_{i-1}}^t f(y(\tau), \tau) d\tau \tag{2.7}$$

; $\forall t \in [t_{i-1}, t_i]$ ; $\forall i \in \bar{J} = \{1, 2, ..., J\}$, with $y(t_0)$, so that, by choosing the real $c_i \in [t_{i-1}, t_i]$ as $c_i = t_i$ for $i = 0, 1, ...., J-1$, one gets from (2.2) in the proof of Lemma 2.1 into (2.3):

$$y(t) = y(t_0) + \sum_{j=0}^{i-2} \int_0^{h_j} \left[ \sum_{k=0}^n \frac{f^{(k)}(y(t_j), t_j)}{k!} (y(\tau + t_j) - y(t_j))^k + \frac{1}{n!} \int_0^\tau (y(\sigma + t_j) - y(t_j))^n f^{(n+1)}(y(\sigma + t_j), \sigma + t_j) d\sigma \right] d\tau$$

$$+ \int_0^{t - t_{i-1}} \left[ \sum_{k=0}^n \frac{f^{(k)}(y(t_{i-1}), t_{i-1})}{k!} (y(\tau + t_{i-1}) - y(t_{i-1}))^k + \frac{1}{n!} \int_0^\tau (y(\sigma + t_{i-1}) - y(t_{i-1}))^n f^{(n+1)}(y(\sigma + t_{i-1}), \sigma + t_{i-1}) d\sigma \right] d\tau$$

(2.8)

; $\forall t \in [t_{i-1}, t_i]$ ; $\forall i \in \bar{J} = \{1, 2, ..., J\}$. Note that, since $f \in C^{(n+1)}(\mathbf{R}^n \times (t_0, t_J); \mathbf{R}^n)$, then $f \in C^{(\ell+1)}(\mathbf{R}^n \times (t_0, t_J); \mathbf{R}^n)$ for any nonnegative integer $\ell \leq n$. Thus, we obtain the result from a similar expression to (2.8) by replacing $n$ by $\ell (\leq n)$ as truncating the Taylor series expansion by its $(\ell + 1)$-th term.

□

A consequence of Theorem 2.2 by using the same technique of the solution construction is the following:

**Corollary 2.3**. Consider the nonlinear ordinary differential equation (2.4) with initial condition $y(t_0)$ on the real interval $\mathbf{R}^n \times \mathbf{R}_{0+}$, with initial conditions $y^{(j)}(t_0)$ for $j = 0, 1, ..., n-1$, such that $f \in C^{(n+1)}(\mathbf{R}^n \times (t_0, t_J); \mathbf{R}^n)$ is Lipschitz-continuous in $[y(t_0) - \theta_0, y(t_0) + \theta_0] \times [t_0, t_J]$ for some $\theta_0 \in \mathbf{R}^n$, and consider also its $\ell - th$ order truncation

$$\dot{x}(t) = \sum_{k=0}^\ell \frac{f^{(k)}(x(t_i), t_i)}{k!} (x(t) - x(t_i))^k \; ; \; x(t_0) = x_0 \tag{2.9}$$

such that $f^{(k)}(y(t),t)$ are bounded in $[y(t_0)-\theta, y(t_0)+\theta] \times [t_0, t_J]$ for $k = 0,1,..., \ell+1$ for some nonnegative integer $\ell \leq n$ and some $\theta \in \overline{B}(\theta_0, R)$ where $\overline{B}(\theta, R) = \{z \in \mathbf{R}^n : \|z - \theta_0\| \leq R\}$ for some positive real $R$ with $x^{(j)}(t_0) = y^{(j)}(t_0)$ for $j = 0,1,..., \ell+1$. □

Since $f^{(k)}(y(t),t)$ are bounded in $[y(t_0)-\theta, y(t_0)+\theta] \times [t_0, \bar{t}_J]$ for $k = 0,1,..., \ell-1$ then, the right-hand-side of (2.9) is Lipschitz-continuous in $[y(t_0)-\theta, y(t_0)+\theta] \times [t_0, t_J] \subseteq [y(t_0)-\theta, y(t_0)+\theta] \times [t_0, t_J]$. Therefore, the unique solution of the truncated differential equation (2.9) in $[a,b]$ is

$$x(t) = x(t_0) + \sum_{k=0}^{\ell} \left\{ \sum_{j=0}^{i-2} \int_0^{h_j} \left[ \frac{f^{(k)}(x(t_j), t_j)}{k!} (x(\tau + t_j) - x(t_j))^k \right] d\tau \right.$$

$$\left. + \int_0^{t-t_{i-1}} \left[ \frac{f^{(k)}(x(t_{i-1}), t_{i-1})}{k!} (x(\sigma + t_{i-1}) - x(t_{i-1}))^k \right] d\tau \right\} \quad (2.10)$$

; $\forall t \in [t_{i-1}, t_i]$ ; $\forall i \in \bar{J} = \{1, 2, ..., J\}$, $\forall \ell (\in \mathbf{Z}_{0+}) \leq n$, where $t_i \in [t_0, t_J]$ are arbitrary strictly ordered points such that $t_0 < t_1 < t_2 < .... < t_{J-1} < t_J$ with $h_i = t_{i+1} - t_i$ for $i = 0,1,...., J-1$. The error in-between the exact solution of (2.9) and that of its truncated form (2.4) is:

$$e(t) = y(t) - x(t)$$

$$= \sum_{j=0}^{i-2} \sum_{k=0}^{\ell} \int_0^{h_j} \left( \frac{f^{(k)}(y(t_j), t_j)}{k!} (y(\tau + t_j) - y(t_j))^k - \frac{f^{(k)}(x(t_j), t_j)}{k!} (x(\tau + t_j) - x(t_j))^k \right) d\tau$$

$$+ \frac{1}{\ell!} \sum_{j=0}^{i-2} \int_0^{h_j} \int_0^{\tau} (y(\sigma + t_j) - y(t_j))^{\ell} f^{(\ell+1)}(y(\sigma + t_j), \sigma + t_j) d\sigma d\tau$$

$$+ \sum_{k=0}^{\ell} \int_0^{t-t_{i-1}} \left( \frac{f^{(k)}(y(t_{i-1}), t_{i-1})}{k!} (y(\tau + t_{i-1}) - y(t_{i-1}))^k - \frac{f^{(k)}(x(t_{i-1}), t_{i-1})}{k!} (x(\tau + t_{i-1}) - x(t_{i-1}))^k \right) d\tau$$

$$+ \frac{1}{\ell!} \int_0^{t-t_{i-1}} \int_0^{\tau} (y(\sigma + t_{i-1}) - y(t_{i-1}))^{\ell} f^{(\ell+1)}(y(\sigma + t_{i-1}), \sigma + t_{i-1}) d\sigma d\tau \quad (2.11)$$

; $\forall t \in [t_{i-1}, t_i]$ ; $\forall i \in \bar{J} = \{1, 2, ..., J\}$, $\forall \ell (\in \mathbf{Z}_{0+}) \leq n$.

*Proof*: Property (i) follows directly Theorem 2.2 applied to the truncated differential equation (2.9) leading to the solution (2.10) in $[t_0, t_J]$. Property (ii) follows from (2.5) and (2.10). □

Now, a preparatory result follows to be then used to guarantee sufficiency-type errors results in-between the true and the approximate solutions in the interval $[a,b]$:

**Lemma 2.4**. Assume that the following hypothesis holds:

A1) $f(y(t), t)$ and its first $(\ell+1)$ derivatives are uniformly bounded from above on a bounded subset of their existence domain with the specific boundedness constraint:

$$\sup_{y(t)\in[y(t_0)-\theta,y(t_J)+\theta],\, t\in[t_0,t_J]} \|f(y(t),t)\| \le K \sup_{y(t)\in[y(t_0)-\theta,y(t_J)+\theta],\, t\in[t_0,t_J]} \|y(t)\| + K_1 \quad (2.12a)$$

$$\sup_{y(t)\in[y(t_0)-\theta,y(t_J)+\theta],\, t\in[t_0,t_J]} \|f^{(j)}(y(t),t)\| \le K \sup_{y(t)\in[y(t_0)-\theta,y(t_J)+\theta],\, t\in[t_0,t_J]} \|f^{(j-1)}(y(t),t)\| + K_1 \quad (2.12b)$$

; for $j = 0, 1, \ldots, \ell+1$ and some $K, K_1 \in \mathbf{R}_{0+}$ with $K < 1$ if $K_1 \in \mathbf{R}_+$. Then, the following properties hold:

**(i)** Assume that the inter-sample intervals $h_i = t_{i+1} - t_i$ for $i = 0, 1, \ldots, J-1$ fulfil the constraint:

$$h_i \le h \le \min\left(a_i,\, \frac{1-\rho_x/2}{\Lambda_x(1-\rho_x/2)^\ell},\, \frac{\rho_x}{2\Lambda_x(1-\rho_x/2)^{\ell+1}(J-1+\rho_x/2)}\right) \quad (2.13)$$

for $i = 0, 1, \ldots, J-1$ and any given real constant $\rho_x \in (0, 2)$, where

$$\Lambda_x = \sum_{k=0}^{\ell} \frac{K^k}{k!}\left(\sup_{y(t)\in[y(t_0)-\theta,y(t_J)+\theta],\, t\in[t_0,t_J]} \|f(y(t),t)\|\right) + \frac{K_1(1-K^{\ell+1})}{1-K} \quad (2.14)$$

$$a_i := \min\, arg\,(t(\in \mathbf{R}_+) > t_i : \|x(t)-x(t_i)\| \le \rho_x/2);\ \forall t \in [t_i, t_{i+1}),\ i = 0, 1, \ldots, J-1 \quad (2.15)$$

Then, the approximated solution fulfils $\sup_{t\in[t_0,t_J]} \|x(t)-x(t_0)\| \le \rho_x/2$ provided that

$$t_1 = \min\, arg\,(t(\in \mathbf{R}_+) > t_0 : \|x(t)-x(t_0)\| \le \rho_x/2)$$

**(ii)** Assume that the inter-sample intervals $h_i = t_{i+1} - t_i$ for $i = 0, 1, \ldots, J-1$ fulfil the constraint:

$$h_i \le h \le \min\left(b_i,\, \frac{1-\rho/2}{\Lambda(1-\rho/2)^\ell},\, \frac{\rho}{2\Lambda(1-\rho/2)^{\ell+1}(J-1+\rho/2)}\right) \quad (2.16)$$

for $i = 0, 1, \ldots, J-1$ and any given real constant $\rho \in (0, 1)$, where

$$\Lambda = \sum_{k=0}^{\ell+1} \frac{K^k}{k!}\left[\left(\sup_{y(t)\in[y(t_0)-\theta,y(t_J)+\theta],\, t\in[t_0,t_J]} \|f(y(t),t)\|\right) + \frac{K_1(1-K^k)}{1-K}\right] \quad (2.17)$$

$$b_i := \min\, arg\,(t > t_i : \|y(t)-y(t_i)\| \le \rho/2);\ \forall t \in [t_i, t_{i+1}),\ i = 0, 1, \ldots, J-1 \quad (2.18)$$

Then, the true solution fulfils $\sup_{t\in[t_0,t_J]} \|y(t)-y(t_0)\| \le \rho/2$ provided that

$$t_1 = \min\, arg\,(t(\in \mathbf{R}_+) > t_0 : \|y(t)-y(t_0)\| \le \rho/2).$$

**(iii)** If $\rho_x = \rho \in (0, 1)$ and, furthermore,

$$h_i \le h \le \min\left(c_i,\, \frac{1-\rho/2}{\Lambda(1-\rho/2)^\ell},\, \frac{\rho}{2\Lambda(1-\rho/2)^{\ell+1}(J-1+\rho/2)}\right);\ i = 0, 1, \ldots, J-1 \quad (2.19)$$

$$c_i := \min\, arg\,(t > t_i : \max(\|y(t)-y(t_i)\|, \|x(t)-x(t_i)\|) \le \rho/2);\ \forall t \in [t_i, t_{i+1}),\ i = 0, 1, \ldots, J-1 \quad (2.20)$$

$$t_1 = \min \arg \left( t(\in \mathbf{R}_+) > t_0 : \max \left( \| y(t) - y(t_0) \|, \| x(t) - x(t_0) \| \right) \leq \rho / 2 \right) \tag{2.21}$$

Then, the true, approximated and error solution fulfil:

$$\sup_{t \in [t_0, t_J]} \| y(t) - y(t_0) \| \leq \rho / 2, \quad \sup_{t \in [t_0, t_J]} \| x(t) - x(t_0) \| \leq \rho / 2 \tag{2.22}$$

and the error in-between them $e(t) = y(t) - x(t)$ fulfils:

$$\sup_{t \in [t_0, t_J]} \| e(t) - e(t_0) \| \leq \rho \tag{2.23}$$

*Proof*: Proceeding recursively, one gets from Assumption A1 that:

$$\sup_{y(t) \in [y(t_0) - \theta, y(t_J) + \theta], t \in [t_0, t_J]} \| f^{(j)}(y(t), t) \| \leq K \sup_{y(t) \in [y(t_0) - \theta, y(t_J) + \theta], t \in [t_0, t_J]} \| f^{(j-1)}(y(t), t) \| + K_1$$

$$\leq K^2 \sup_{y(t) \in [y(t_0) - \theta, y(t_J) + \theta], t \in [t_0, t_J]} \| f^{(j-1)}(y(t), t) \| + K_1(1 + K)$$

$$\ldots\ldots\ldots\ldots\ldots\ldots\ldots\ldots\ldots\ldots\ldots$$

$$\leq K^k F_0 + K_1 \left( \sum_{i=0}^{k-1} K^i \right) \leq K^k F_0 + \frac{K_1 (1 - K^k)}{1 - K} \leq F_0 + \frac{K_1}{1 - K} \tag{2.24}$$

if $K < 1$ and $K_1 \neq 0$, and

$$\sup_{y(t) \in [y(t_0) - \theta, y(t_J) + \theta], t \in [t_0, t_J]} \| f^{(k)}(y(t), t) \| \leq K^k F_0 \tag{2.25}$$

If $K_1 = 0$, where

$$F_0 = \sup_{y(t) \in [y(t_0) - \theta, y(t_J) + \theta], t \in [t_0, t_J]} \| f^{(0)}(y(t), t) \| = \sup_{y(t) \in [y(t_0) - \theta, y(t_J) + \theta], t \in [t_0, t_J]} \| f(y(t), t) \| < +\infty \tag{2.26}$$

*Case a*: if $K < 1$ and $K_1 \neq 0$ proceed by complete induction by assuming that $\sup_{t \in [t_0, t_i]} \| x(t) - x(t_0) \| \leq \rho_x / 2$ since the condition $\left( t(\in \mathbf{R}_+) > t_i : \| x(t) - x(t_i) \| \leq \rho_x / 2 \right)$ guarantees that $\sup_{t \in [t_0, t_1]} \| x(t) - x(t_0) \| \leq \rho_x / 2$. Thus, one gets from (2.10) that

$$\| x(t) - x(t_i) \| \leq \sum_{k=0}^{\ell} \sum_{j=0}^{i} \frac{h_j}{k!} \left( K^k F_0 + \frac{K_1 (1 - K^k)}{1 - K} \right) (\rho_x / 2)^k \tag{2.27}$$

$$= \sum_{k=0}^{\ell} \sum_{j=0}^{i-1} \frac{h_j}{k!} \left( K^k F_0 + \frac{K_1 (1 - K^k)}{1 - K} \right) (\rho_x / 2)^k + \sum_{k=0}^{\ell} \frac{h_i}{k!} \left( K^k F_0 + \frac{K_1 (1 - K^k)}{1 - K} \right) \| x(t) - x(t_i) \|^k$$

$$\leq \sum_{k=0}^{\ell} \sum_{j=0}^{i-1} \frac{h_j}{k!} \left( K^k F_0 + \frac{K_1 (1 - K^k)}{1 - K} \right) (\rho_x / 2)^k + \sum_{k=0}^{\ell} \frac{h_i}{k!} \left( K^k F_0 + \frac{K_1 (1 - K^k)}{1 - K} \right) \| x(t) - x(t_i) \|^k$$

$$= \Lambda_x \left( \sum_{k=0}^{\ell} (\rho_x / 2)^k \right) \left( \sum_{j=0}^{i-1} h_j \right) + h_i \Lambda \left[ \sum_{k=0}^{\ell} \left( \| x(t) - x(t_i) \|^k \right) \right]$$

$$\leq i h \Lambda_x \left( \sum_{k=0}^{\ell} (\rho_x / 2)^k \right) + h_i \Lambda_x \left( \sum_{k=0}^{\ell} \left( \| x(t) - x(t_i) \|^{k-1} \right) \right) \| x(t) - x(t_i) \|; \ \forall t \in [t_i, t_{i+1}) \tag{2.28}$$

where $\Lambda_x = \sum_{k=0}^{\ell} \frac{1}{k!}\left( K^k F_0 + \frac{K_1(1-K^k)}{1-K} \right)$ and $h \geq \max_{0 \leq i \leq J-1} h_i$, with $h_i = t_{i+1} - t_i$, for $i = 0, 1, ..., J-1$, so that

$$\left[1 - h_i \Lambda_x \left(\sum_{k=0}^{\ell}\left(\|x(t)-x(t_i)\|^{k-1}\right)\right)\right]\|x(t)-x(t_i)\| \leq ih\Lambda_x\left(\sum_{k=0}^{\ell}(\rho_x/2)^k\right) = \frac{ih\Lambda_x\left(1-(\rho_x/2)^{\ell+1}\right)}{1-\rho_x/2} \quad (2.29)$$

or

$$\|x(t)-x(t_i)\| \leq \frac{ih\Lambda_x\left(1-(\rho_x/2)^{\ell+1}\right)}{\left[1 - h_i\Lambda_x\left(\sum_{k=0}^{\ell}\left(\|x(t)-x(t_i)\|^{k-1}\right)\right)\right](1-\rho_x/2)} \leq \frac{ih\Lambda_x\left(1-(\rho_x/2)^{\ell+1}\right)}{\left[1 - \frac{h_i\Lambda_x\left(1-(\rho_x/2)^{\ell}\right)}{1-\rho_x/2}\right](1-\rho_x/2)} \leq \frac{\rho_x}{2} \quad (2.30)$$

provided that $0 < \rho_x < 2$, and $1 > h_i \Lambda_x \left(\sum_{k=0}^{\ell}\left(\|x(t)-x(t_i)\|^{k-1}\right)\right)$ which is guaranteed if $h_i < \frac{1-\rho_x/2}{\Lambda_x(1-\rho_x/2)^{\ell}}$ holds with $a_i$ for $i = 0, 1, ..., J-1$ being defined in (2.15), provided that

$\|x(t)-x(t_j)\| \leq \rho_x/2$; $\forall t \in [t_j, t_{j+1})$ for $j(\leq i) = 0, 1, ..., i-1$ and, then (2.30) and $h_j < \frac{1-\rho_x/2}{\Lambda_x(1-\rho_x/2)^{\ell}}$ for

$j = 0, 1, ..., i-1$ are jointly guaranteed for the given $i = 0, 1, ..., J-1$ if

$$h_i < \min\left( \frac{1-\rho_x/2}{\Lambda_x(1-\rho_x/2)^{\ell}}, \frac{\rho_x}{2\Lambda_x(1-\rho_x/2)^{\ell+1}(J-1+\rho_x/2)} \right) \quad (2.31)$$

provided that $\|x(t)-x(t_j)\| \leq \rho_x/2$ for $t \in [t_j, t_{j+1})$ for $j = 0, 1, ..., i-1$, the last condition being identical to

$$t_{i+1} \leq a_i := \min arg\left(t > t_i : \|x(t)-x(t_i)\| \leq \rho_x/2\right) \quad (2.32)$$

The above two conditions (2.31) - (2.32) jointly reduce to (2.13). Then, one gets from complete induction from (2.28), if (2.13) holds, that:

$\sup_{t \in [t_0, t_i]} \|x(t)-x(t_i)\| \leq \rho_x/2 \Rightarrow \sup_{t \in [t_0, t_{i+1})} \|x(t)-x(t_i)\| \leq \rho_x/2$ and, one gets also by continuity extension,

$\sup_{t \in [t_0, t]} \|x(t)-x(t_0)\| \leq \rho_x/2$; $\forall t \in [t_0, t_J]$. Hence, we have got the result for Case a.

*Case b*: If $K_1 = 0$ then

$$\|x(t)-x(t_i)\| \leq ih\Lambda_{x0}\left(\sum_{k=0}^{\ell}(\rho_x/2)^k\right) + h_i\Lambda_{x0}\left(\sum_{k=0}^{\ell}\left(\|x(t)-x(t_i)\|^{k-1}\right)\right)\|x(t)-x(t_i)\| \quad (2.33)$$

; $\forall t \in [t_i, t_{i+1})$, where $\Lambda_{x0} = \sum_{k=0}^{\ell} \frac{K^k F_0}{k!} \leq \Lambda_x$, so that $(1-h\Lambda_{x0})\|x(t)-x(t_i)\| \leq ih\Lambda_{x0} E$ for $i = 0, 1, ..., J$ and thus, one gets:

$$\|x(t)-x(t_i)\| \leq \frac{ih\Lambda_{x0}E}{1-h\Lambda_{x0}} \leq \rho_x/2 \; ; \; \forall t \in [t_i, t_{i+1}) \text{ for } i=0,1,\ldots,J$$ and one gets from complete induction the same conclusion $\sup_{t \in [t_0, t]} \|x(t)-x(t_0)\| \leq \rho_x/2 \; ; \; \forall t \in [t_0, t_J]$ as in Case a provided that (2.13) holds.

Then, (2.13) guarantees Property (i) for both Case a and Case b. Then, Property (i) has been proven. Property (ii) is proven "mutatis-mutandis" by noting that $\Lambda \geq \Lambda_x$ from (2.14) and (2.17) and noting also that $a_i$ in (2.15) is replaced with $b_i$ in (2.18) so that the admissible inter-sample interval satisfying the constraint (2.13) is replaced by such an interval satisfying the constraint (2.16). Finally, Property (iii) follows from Properties (i)-(ii) by equalizing $\rho_x$ and $\rho$ to take a maximum value being less than $1/2$. □

The following comments address the fact that it is not necessary to deal with the solution of the true differential equation (2.4) to calculate in Lemma 2.4.

**Remark 2.5**: Note that one gets by direct integration from (2.4) that:

$$\|y(t)\| \leq \sup_{t_0 \leq \tau \leq t_J} \|y(\tau)\| \leq \|y(t_0)\| + (t_J - t_0)K \sup_{t_0 \leq \tau \leq t_J} \|y(\tau)\| \tag{2.34}$$

leading to

$$\sup_{t_0 \leq t \leq t_J} \|y(t)\| \leq \frac{K(t_J - t_0)}{1 - K(t_J - t_0)} \|y(t_0)\| \text{ if } \sum_{i=0}^{J-1} h_i < 1/K.$$

Thus, (2.22)-(2.23) might be guaranteed with $\sup_{t_0 \leq t \leq t_J} \|y(t)-y(t_0)\| \leq \frac{1}{1-K(t_J-t_0)} \|y(t_0)\| \leq \frac{\rho}{2}$ if $\|y(t_0)\| \leq \frac{\rho}{2}(1-K(t_J-t_0)) < \frac{\rho}{2}$. Thus, there is no need to compute the solution of the true differential equation (2.4) and $\sup_{t_i \leq t < t_{i+1}} \|y(t)-y(t_0)\| \leq \frac{\rho}{2}$ for $i=0,1,\ldots,J-1$ in (2.18) and (2.20) if

$$\|y(t_0)\| \leq \frac{\rho}{2}(1-K(t_J-t_0)). \qquad \square$$

One gets directly from Lemma 2.4 the subsequent result:

**Theorem 2.6**. Assume that the conditions (2.12) and (2.19)-(2.21) of Lemma 2.4 (iii). Then

$$max\left(\max_{0 \leq i \leq J-1} \|e(t_{i+1})-e(t_i)\|, \|e(t)-e(t_0)\|\right) \leq \rho < 1 \; ; \; \max_{t \in [t_0, t_J]} \|e(t)\| \leq \|e(t_0)\| + \rho \tag{2.35}$$

$; \forall t \in [t_0, t_J]$, and

$$\max_{0 \leq i \leq J-1} \|e(t)-e(t_i)\| \leq \rho \; ; \; \forall t \in [t_i, t_{i+1}] \text{ for } i=0,1,\ldots,J-1 \tag{2.36}$$

for $i=0,1,\ldots,J-1$. □

## 3. Orbits of the exact solution, pseudo-orbits of the approximated solution and the shadowing property

Now, consider a perturbed solution (2.10) of the approximated differential equation (2.9) associated with a perturbation $x(t_i) = x(t_i^-) + g(t_i)$ with $\{g(t_i)\} \subset \mathbf{R}^n$ at $t = t_i$ fulfilling $\|g(t_i)\| \leq \bar{g}_i \leq \bar{g}$ for some given $\bar{g} \in \mathbf{R}$, $\forall i \in \bar{J}$. Note that a perturbation at the initial $t = t_0$ is considered by choosing $x(t_0) = y(t_0) + g_0$ for some nonzero $g_0 = g(t_0) \in \mathbf{R}$. The perturbed solution can be generated, in particular, from impulsive controls of amplitudes $g(t_i)$ at the sequence $\{t_i : i \in \bar{J}\}$. The exact and approximate solutions (2.5) and (2.10) in $[t_0, t_J]$, provided that they exist, are:

$$y(t) = y(t_i) + \sum_{k=0}^{\ell} \int_0^{t-t_i} \frac{f^{(k)}(y(t_i), t_i)}{k!} (y(\tau + t_i) - y(t_i))^k \, d\tau + \frac{1}{\ell!} \int_0^{t-t_i} \int_0^{\tau} (y(\sigma + t_j) - y(t_j))^{\ell} f^{(\ell+1)}(y(\sigma + t_j), \sigma + t_j) \, d\sigma \, d\tau$$

$; \forall t \in [t_i, t_{i+1}], \ i = 0, 1, \ldots, J-1$ \hfill (3.3)

$$x(t) = x(t_i) + \sum_{k=0}^{\ell} \int_0^{t-t_i} \frac{f^{(k)}(x(t_i), t_i)}{k!} (x(\tau + t_i) - x(t_i))^k \, d\tau + U(t - t_{i+1}) g(t_{i+1})$$

$; \forall t \in [t_i, t_{i+1}], \ i = 0, 1, \ldots, J-1$ \hfill (3.4)

where $U(t)$ is the Heaviside function. The error between the exact and perturbed approximated solutions becomes:

$$e(t) = e(t_i) + \sum_{k=0}^{\ell} \int_0^{t-t_i} \left( \frac{f^{(k)}(y(t_i), t_i)}{k!} (y(\tau + t_i) - y(t_i))^k - \frac{f^{(k)}(x(t_i), t_i)}{k!} (x(\tau + t_i) - x(t_i))^k \right) d\tau$$

$$+ \frac{1}{\ell!} \int_0^{t-t_i} \int_0^{\tau} (y(\sigma + t_j) - y(t_j))^{\ell} f^{(\ell+1)}(y(\sigma + t_j), \sigma + t_j) \, d\sigma \, d\tau - U(t - t_{i+1}) g(t_{i+1})$$ \hfill (3.5)

$; \forall t \in [t_i, t_{i+1}]; \ i = 0, 1, \ldots, J-1$. Now, one gets from (2.22) - (2.23) of Lemma 2.4:

$$\|e(t) - e(t_i)\| \leq \sum_{k=0}^{\ell} \frac{2^{k+1}}{k!} (t - t_i) M_{ik} \left( \frac{\rho}{2} \right)^k + \frac{1}{\ell!} (t - t_i)^2 2^{\ell} M_{i,\ell+1} \left( \frac{\rho}{2} \right)^{\ell} + \bar{g}_i \ ; \ \forall t \in [t_i, t_{i+1}]$$ \hfill (3.6)

where

$$M_{ik} = \sup_{y(t) \in [y(t_i) - \theta, y(t_i) + \theta], t \in [t_i, t_{i+1}]} \|f^{(j)}(y(t), t)\| \leq K^k F_0 + \frac{K_1(1 - K^k)}{1 - K}$$

$$\leq \frac{K^{k+1} \rho}{2} + \frac{K_1(1 - K^k)}{1 - K} \leq \frac{\rho}{2} + \frac{K_1}{1 - K}$$ \hfill (3.7)

from (2.24) and, one gets after using the triangle inequality:

$$\|e(t) - e(t_i)\| \leq \sum_{j=i}^{m} \left( \sum_{k=0}^{\ell} \frac{2^{k+1}}{k!} (t - t_j) \left[ \frac{K^{k+1} \rho}{2} + \frac{K_1(1 - K^k)}{1 - K} \right] \left( \frac{\rho}{2} \right)^k \right.$$

$$\left. + \frac{1}{\ell!} (t - t_j)^2 2^{\ell} \left[ \frac{K^{\ell+2} \rho}{2} + \frac{K_1(1 - K^{\ell+1})}{1 - K} \right] \left( \frac{\rho}{2} \right)^{\ell} + \bar{g}_j \right)$$ \hfill (3.8)

$; \forall t \in [t_{i+m}, t_{i+m+1}]$ for $m = 0, 1, \ldots, J - i - 1$; $i = 0, 1, \ldots, J - 1$. One obtains easily from (3.8), either by complete induction or via recursive calculation, that:

$$\|e(t)-e(t_i)\| \le \sum_{j=0}^{J-1}\left(\sum_{k=0}^{\ell}\frac{2^{k+1}}{k!}(t-t_j)\left[\frac{K^{k+1}\rho}{2}+\frac{K_1(1-K^k)}{1-K}\right]\left(\frac{\rho}{2}\right)^k\right.$$

$$\left.+\frac{1}{\ell!}(t-t_j)^2 2^{\ell}\left[\frac{K^{\ell+2}\rho}{2}+\frac{K_1(1-K^{\ell+1})}{1-K}\right]\left(\frac{\rho}{2}\right)^{\ell}+\overline{g}_j\right) \quad (3.9a)$$

$$\le Jh\left(\sum_{k=0}^{\ell}\frac{2^{k+1}}{k!}\left[\frac{K^{k+1}\rho}{2}+\frac{K_1(1-K^k)}{1-K}\right]\left(\frac{\rho}{2}\right)^k+\frac{Jh}{\ell!}2^{\ell}\left[\frac{K^{\ell+2}\rho}{2}+\frac{K_1(1-K^{\ell+1})}{1-K}\right]\left(\frac{\rho}{2}\right)^{\ell}\right)+J\overline{g} \quad (3.9b)$$

$$\|e(t)\| \le \|e(t_0)\|+\sum_{j=0}^{J-1}(t-t_j)\left(\sum_{k=0}^{\ell}\frac{2^{k+1}}{k!}\left[\frac{K^{k+1}\rho}{2}+\frac{K_1(1-K^k)}{1-K}\right]\left(\frac{\rho}{2}\right)^k\right.$$

$$\left.+\frac{1}{\ell!}(t-t_j)2^{\ell}\left[\frac{K^{\ell+2}\rho}{2}+\frac{K_1(1-K^{\ell+1})}{1-K}\right]\left(\frac{\rho}{2}\right)^{\ell}+\overline{g}_j\right) \quad (3.10a)$$

$$\le \|e(t_0)\|+Jh\left(\sum_{k=0}^{\ell}\frac{2^{k+1}}{k!}\left[\frac{K^{k+1}\rho}{2}+\frac{K_1(1-K^k)}{1-K}\right]\left(\frac{\rho}{2}\right)^k+\frac{Jh}{\ell!}2^{\ell}\left[\frac{K^{\ell+2}\rho}{2}+\frac{K_1(1-K^{\ell+1})}{1-K}\right]\left(\frac{\rho}{2}\right)^{\ell}\right)+J\overline{g}$$

(3.10b)

; $\forall t\in\left[t_0,t_0+\sum_{i=0}^{J-1}h_i\right]\left(\subseteq\left[t_0,t_0+Jh\right]\right)$ with $h_i=t_{i+1}-t_i$ for $i=0,1,\ldots,J-1$ and any given nonnegative integer $\ell\le n$. The following result follows directly from the above equations and Theorem 2.6:

**Theorem 3.1**. Consider an approximated perturbed solution (3.4) under a forcing perturbation sequence $\{g(t_i)\}\subset \boldsymbol{R}^n$ at $t=t_i$ fulfilling $\|g(t_i)\|\le \overline{g}_i\le \overline{g}\ge\|e(t_0)\|$ for $i=1,2,\ldots$ and some $\overline{g}\in \boldsymbol{R}_+$. Then, there are numbers $h\in\boldsymbol{R}_+$, $J=J(h)\in\boldsymbol{Z}_+$, $\varepsilon_1=\varepsilon_1(h,\overline{g})\in\boldsymbol{R}_+$ and $\varepsilon=\varepsilon(\varepsilon_1,\|e(t_0)\|)$ such that

$$max\left(\max_{0\le i\le J-1}\|e(t_{i+1})-e(t_i)\|, \max_{t\in\boldsymbol{R}_{0+}}\|e(t)-e(t_0)\|\right)\le \varepsilon_1 \; ; \; \max_{t\in\boldsymbol{R}_{0+}}\|e(t)\|\le \varepsilon \quad (3.11)$$

on $[t_0,t_J]$ for any strictly ordered sequence of $(J+1)$ nonnegative real numbers $\{t_i:i=0,1,\ldots,J\}$, subject to the constraints:

$$t_J=t_0+\sum_{i=0}^{J-1}h_i \quad , \qquad h_i=t_{i+1}-t_i\le h \; ; \; i=0,1,\ldots,J-1 \quad (3.12)$$

satisfying the constraints (2.19)-(2.21) of Lemma 2.4 subject to (2.16).

*Proof*: Note that fixing $\sum_{i=0}^{J-1}h_i=t_J-t_0\le Jh$, with $h=\max_{0\le i\le J-1}(t_{i+1}-t_i)$, and the use of (3.8)-(3.9) leads to

$$\|e(t)-e(t_i)\| \le \sum_{i=0}^{J-1}\left(h_i\left(\sum_{k=0}^{\ell}\frac{2^{k+1}}{k!}\left[\frac{K^{k+1}\rho}{2}+\frac{K_1(1-K^k)}{1-K}\right]\left(\frac{\rho}{2}\right)^k+\frac{1}{\ell!}h2^{\ell}\left[\frac{K^{\ell+2}\rho}{2}+\frac{K_1(1-K^{\ell+1})}{1-K}\right]\left(\frac{\rho}{2}\right)^{\ell}\right)+\overline{g}_i\right)$$

$$\le \varepsilon_1=\rho+\sum_{i=0}^{J-1}\overline{g}_i \quad (3.13a)$$

$$\|e(t)-e(t_i)\| \le Jh\left(\sum_{k=0}^{\ell}\frac{2^{k+1}}{k!}\left[\frac{K^{k+1}\rho}{2}+\frac{K_1(1-K^k)}{1-K}\right]\left(\frac{\rho}{2}\right)^k+\frac{1}{\ell!}h2^{\ell}\left[\frac{K^{\ell+2}\rho}{2}+\frac{K_1(1-K^{\ell+1})}{1-K}\right]\left(\frac{\rho}{2}\right)^{\ell}\right)$$

$$+J\bar{g} \leq \bar{\varepsilon}_1 = \bar{\rho} + J\bar{g} \tag{3.13b}$$

; $\forall t \in [t_0, t_J]$; $i = 0, 1, ..., J-1$ applying Lemma 2.4 and Theorem 2.6 for any given prefixed $\rho \in \mathbf{R}_+$.

The result then follows since $\bar{g}_0 = \|e(t_0)\| \leq \bar{g}$ and either

$$\varepsilon = \varepsilon_1 + \|e(t_0)\|, \quad \varepsilon_1 = \rho + \sum_{i=0}^{J-1}\bar{g}_i,$$

$$\sum_{i=0}^{J-1} h_i \left( \sum_{k=0}^{\ell} \frac{2^{k+1}}{k!} \left[ \frac{K^{k+1}\rho}{2} + \frac{K_1(1-K^k)}{1-K} \right] \left(\frac{\rho}{2}\right)^k + \frac{h}{\ell!} 2^\ell \left[ \frac{K^{\ell+2}\rho}{2} + \frac{K_1(1-K^{\ell+1})}{1-K} \right] \left(\frac{\rho}{2}\right)^\ell \right) \leq \rho \tag{3.14}$$

or

$$\varepsilon = \bar{\varepsilon}_1 + \|e(t_0)\|, \quad \bar{\varepsilon}_1 = \bar{\rho} + J\bar{g},$$

$$h \left( \sum_{k=0}^{\ell} \frac{2^{k+1}}{k!} \left[ \frac{K^{k+1}\rho}{2} + \frac{K_1(1-K^k)}{1-K} \right] \left(\frac{\rho}{2}\right)^k + \frac{1}{\ell!} h 2^\ell \left[ \frac{K^{\ell+2}\rho}{2} + \frac{K_1(1-K^{\ell+1})}{1-K} \right] \left(\frac{\rho}{2}\right)^\ell \right) \leq \bar{\rho} \tag{3.15}$$

and the result has been proven. □

The following result extends Theorem 3.1 with results of Theorem 2.6 for the case when both the exact and approximated differential equations are subject to a piecewise- continuous bounded disturbance which might be dependent on the solution and also can have finite step discontinuities in the sequence $\{t_i : i = 0, 1, ..., J\}$:

**Theorem 3.2**. Consider the forced versions of the differential equations (2.4) and (2.9):

$$\dot{y}(t) = f(y(t), t) + g(\tau, y(\tau)), \quad y(t_0) = y_0 \tag{3.16}$$

$$\dot{x}(t) = \sum_{k=0}^{\ell} \frac{f^{(k)}(x(t_i), t_i)}{k!} (x(t) - x(t_i))^k + g(\tau, x(\tau)); \quad x(t_0) = x_0 \tag{3.17}$$

under the additive forcing perturbation $g \in C^{(n+1)}(\mathbf{R}^n \times (t_0, t_J); \mathbf{R}^n)$ satisfying Assumption A2 of Lemma 2.4 fulfilling $g(y(t), t) = \lambda(t) y(t)$ and $g(x(t), t) = \lambda(t) x(t) + U(t - t_{i+1}) g_{i+1}$; $\forall t \in [t_i, t_{i+1}]$ with $\|g_{i+1}\| \leq \bar{g}$ for $i = 0, 1, ..., J-1$ and some $\bar{g} \in \mathbf{R}_+$ and $\lambda : [t_0, t_J] \to \mathbf{R}^n$ being a bounded piecewise continuous function. Then, there are numbers $h \in \mathbf{R}_+$, $J = J(h) \in \mathbf{Z}_+$, $\varepsilon_1 = \varepsilon_1(h, \bar{g}) \in \mathbf{R}_+$ and $\varepsilon = \varepsilon(\varepsilon_1, \|e(t_0)\|)$ such that

$$max\left( \max_{0 \leq i \leq J-1} \|e(t_{i+1}) - e(t_i)\|, \max_{t \in \mathbf{R}_{0+}} \|e(t) - e(t_0)\| \right) \leq \varepsilon_1; \quad \max_{t \in \mathbf{R}_{0+}} \|e(t)\| \leq \rho \tag{3.18}$$

on $[t_0, t_J]$ for a strictly ordered finite set of $(J+1)$ nonnegative real numbers $\{t_i : i = 0, 1, ..., J\}$, subject to the constraints $t_J = t_0 + \sum_{i=0}^{J-1} h_i$, $h_i = t_{i+1} - t_i \leq h$; $i = 0, 1, ..., J-1$, the constraints (2.19)-(2.21) subject to (2.16), and either

$$\sum_{i=0}^{J-1} h_i \lambda_i < 1 \tag{3.19a}$$

$$\sum_{i=0}^{J-1} \bar{g}_i < \infty, \quad \bar{g} \geq \frac{\sum_{i=0}^{J-1} h_i \lambda_i}{1 - \sum_{i=0}^{J-1} h_i \lambda_i} \|e(t_0)\| \tag{3.19b}$$

or

$$Jh\lambda < 1 \tag{3.20.a}$$

$$J\bar{g} < \infty \quad , \quad \bar{g} \geq \frac{Jh\lambda}{1-Jh\lambda}\|e(t_0)\| \tag{3.20.b}$$

*Proof*: Fix $\sum_{i=0}^{J-1} h_i = t_J - t_0 \leq Jh$, with $h = \max_{0 \leq i \leq J-1}(t_{i+1} - t_i)$. Eqs. (3.16)- (3.17) have the following solutions:

$$y(t) = y(t_i) + \sum_{k=0}^{\ell} \int_0^{t-t_i} \left( \frac{f^{(k)}(y(t_i), t_i)}{k!} (y(\tau + t_i) - y(t_i))^k + \frac{1}{\ell!} \int_0^{\tau} (y(\sigma + t_j) - y(t_j))^{\ell} f^{(\ell+1)}(y(\sigma + t_j), \sigma + t_j) d\sigma + g(y(\tau), \tau) \right) d\tau$$

(3.21)

$$x(t) = x(t_i) + \int_0^{t-t_i} \left( \sum_{k=0}^{\ell} \frac{f^{(k)}(x(t_i), t_i)}{k!} (x(\tau + t_i) - x(t_i))^k + g(x(\tau), \tau) \right) d\tau + U(t - t_{i+1}) g_{i+1} \tag{3.22}$$

; $\forall t \in [t_i, t_{i+1}]$, $i = 0, 1, ..., J-1$. Note that:

$$g(y(t), t) - g(x(t), t) = \lambda(t)(y(t) - x(t)) = \lambda(t) e(t)$$

Thus, the error between both of them becomes:

$$e(t) = e(t_i) + \sum_{k=0}^{\ell} \int_0^{t-t_i} \left( \frac{f^{(k)}(y(t_i), t_i)}{k!} (y(\tau + t_i) - y(t_i))^k - \frac{f^{(k)}(x(t_i), t_i)}{k!} (x(\tau + t_i) - x(t_i))^k + g(e(\tau), \tau) \right) d\tau$$

$$+ \frac{1}{\ell!} \int_0^{t-t_i} \int_0^{\tau} (y(\sigma + t_j) - y(t_j))^{\ell} f^{(\ell+1)}(y(\sigma + t_j), \sigma + t_j) d\sigma d\tau - U(t - t_{i+1}) g(x(t_{i+1}), t_{i+1}) \tag{3.23}$$

; $\forall t \in [t_i, t_{i+1}]$. Then, (3.7) leads to:

$$\|e(t)\| \leq \|e(t_i)\| + (t - t_i)\left( \sum_{k=0}^{\ell} \frac{2^{k+1}}{k!} \left[ \frac{K^{k+1}\rho}{2} + \frac{K_1(1-K^k)}{1-K} \right] \left(\frac{\rho}{2}\right)^k \right.$$

$$\left. + \frac{1}{\ell!}(t - t_i) 2^{\ell} \left[ \frac{K^{\ell+2}\rho}{2} + \frac{K_1(1-K^{\ell+1})}{1-K} \right] \left(\frac{\rho}{2}\right)^{\ell} + g_{i+1} \right) \tag{3.24}$$

; $\forall t \in [t_i, t_{i+1}]$, $i = 0, 1, ..., J-1$. Then,

$$\sup_{t_i \leq t \leq t_{i+1}} (\|e(t)\|) \leq \|e(t_i)\| + (t - t_i)\left( \sum_{k=0}^{\ell} \frac{2^{k+1}}{k!} \left[ \frac{K^{k+1}\rho}{2} + \frac{K_1(1-K^k)}{1-K} \right] \left(\frac{\rho}{2}\right)^k \right.$$

$$\left. + \frac{1}{\ell!}(t - t_i) 2^{\ell} \left[ \frac{K^{\ell+2}\rho}{2} + \frac{K_1(1-K^{\ell+1})}{1-K} \right] \left(\frac{\rho}{2}\right)^{\ell} + g_{i+1} \right) + h_i \lambda_i \sup_{t_i \leq \tau \leq t_{i+1}} (\|e(\tau)\|) + \bar{g} \tag{3.25}$$

so that, since $1 > h_i \lambda_i$, where $\lambda_i = \max_{t_i \leq \tau \leq t_{i+1}}(\lambda(\tau))$ for $i = 0, 1, ..., J-1$, one gets:

$$\sup_{t_i \leq t \leq t_{i+1}} (\|e(t)\|) \leq \frac{1}{1 - h_i \lambda_i}$$

$$\times \left( \|e(t_i)\| + \left( \sum_{k=0}^{\ell} \frac{2^{k+1}}{k!} h_i \left[ \frac{K^{k+1}\rho}{2} + \frac{K_1(1-K^k)}{1-K} \right] \left(\frac{\rho}{2}\right)^k + \frac{1}{\ell!}(t-t_i)2^{\ell}\left[ \frac{K^{\ell+2}\rho}{2} + \frac{K_1(1-K^{\ell+1})}{1-K} \right]\left(\frac{\rho}{2}\right)^{\ell} + g_{i+1} \right) \right)$$

(3.26)

what implies

$$\left| \sup_{t_i \leq t \leq t_{i+1}} (\|e(t)\| - \|e(t_i)\|) \right| \leq \frac{h_i \lambda_i}{1 - h_i \lambda_i} \|e(t_i)\| + \frac{1}{1 - h_i \lambda_i}$$

$$\times \left( \sum_{k=0}^{\ell} \frac{2^{k+1}}{k!} h_i \left[ \frac{K^{k+1}\rho}{2} + \frac{K_1(1-K^k)}{1-K} \right] \left(\frac{\rho}{2}\right)^k + \frac{1}{\ell!} h_i^2 2^{\ell}\left[ \frac{K^{\ell+2}\rho}{2} + \frac{K_1(1-K^{\ell+1})}{1-K} \right]\left(\frac{\rho}{2}\right)^{\ell} + g_{i+1} \right)$$

(3.27)

If $\sum_{i=0}^{J-1} h_i \lambda_i < 1$, we also get (3.28)-(3.29) below from (3.27) as well as (3.30)-(3.31) if, in addition, $Jh\lambda < 1$:

$$\sup_{t_0 \leq t \leq t_J} (\|e(t)\|) \leq \frac{1}{1 - \sum_{i=0}^{J-1} h_i \lambda_i}$$

$$\times \left( \|e(t_0)\| + \sum_{i=0}^{J-1} h_i \left( \sum_{k=0}^{\ell} \frac{2^{k+1}}{k!} \left[ \frac{K^{k+1}\rho}{2} + \frac{K_1(1-K^k)}{1-K} \right] \left(\frac{\rho}{2}\right)^k + \frac{1}{\ell!} h_i 2^{\ell}\left[ \frac{K^{\ell+2}\rho}{2} + \frac{K_1(1-K^{\ell+1})}{1-K} \right]\left(\frac{\rho}{2}\right)^{\ell} \right) + \sum_{i=0}^{J-1} \overline{g}_i \right)$$

(3.28)

$$\left| \sup_{t_0 \leq t \leq t_J} (\|e(t)\| - \|e(t_0)\|) \right| \leq \frac{\sum_{i=0}^{J-1} h_i \lambda_i}{1 - \sum_{i=0}^{J-1} h_i \lambda_i}$$

$$\times \left( \|e(t_0)\| + \sum_{i=0}^{J-1} h_i \left( \sum_{k=0}^{\ell} \frac{2^{k+1}}{k!} \left[ \frac{K^{k+1}\rho}{2} + \frac{K_1(1-K^k)}{1-K} \right] \left(\frac{\rho}{2}\right)^k + \frac{1}{\ell!} h_i 2^{\ell}\left[ \frac{K^{\ell+2}\rho}{2} + \frac{K_1(1-K^{\ell+1})}{1-K} \right]\left(\frac{\rho}{2}\right)^{\ell} \right) + \sum_{i=0}^{J-1} \overline{g}_i \right)$$

(3.29)

$$\sup_{t_0 \leq t \leq t_J} (\|e(t)\|) \leq \frac{1}{1 - Jh\lambda}$$

$$\times \left( \|e(t_0)\| + Jh \left( \sum_{k=0}^{\ell} \frac{2^{k+1}}{k!} \left[ \frac{K^{k+1}\rho}{2} + \frac{K_1(1-K^k)}{1-K} \right] \left(\frac{\rho}{2}\right)^k + \frac{1}{\ell!} h 2^{\ell}\left[ \frac{K^{\ell+2}\rho}{2} + \frac{K_1(1-K^{\ell+1})}{1-K} \right]\left(\frac{\rho}{2}\right)^{\ell} \right) + J\overline{g} \right)$$

(3.30)

$$\left| \sup_{t_0 \leq t \leq t_J} (\|e(t)\| - \|e(t_0)\|) \right| \leq \frac{Jh\lambda}{1 - Jh\lambda}$$

$$\times Jh \left( \sum_{k=0}^{\ell} \frac{2^{k+1}}{k!} \left[ \frac{K^{k+1}\rho}{2} + \frac{K_1(1-K^k)}{1-K} \right] \left(\frac{\rho}{2}\right)^k + \frac{1}{\ell!} h 2^{\ell}\left[ \frac{K^{\ell+2}\rho}{2} + \frac{K_1(1-K^{\ell+1})}{1-K} \right]\left(\frac{\rho}{2}\right)^{\ell} \right) + J\overline{g}$$

(3.31)

where $\lambda = \max_{0 \leq t \leq t_J}(\lambda(t)) = \max_{0 \leq i \leq J-1} \max_{t_i \leq \tau \leq t_{i+1}}(\lambda(\tau)) = \max_{0 \leq i \leq J-1} \lambda_i$. Property (i) follows from (3.27)-(3.28) by defining $\varepsilon_0$, $\varepsilon_1$ and $\varepsilon$ as in (3.14) since $\overline{g} \geq \max\left( \max_{0 \leq i \leq J-1} \|g_{i+1}\|, \frac{\sum_{i=0}^{J-1} h_i \lambda_i}{1 - \sum_{i=0}^{J-1} h_i \lambda_i} \|e(t_0)\| \right)$ and Property (ii)

follows from (3.29)-(3.30) by defining $\rho$, $\varepsilon_1$ and $\varepsilon$ as in (3.15) since $\bar{g} \geq max\left( \max_{0 \leq i \leq J-1} \|g_{i+1}\|, \frac{Jh\lambda}{1-Jh\lambda} \|e(t_0)\| \right)$. Thus, the result has been proven. □

Now, three definitions are given concerning the so-called pseudo-orbits, as a counterpart to the true sampled trajectory solution, or orbit, of finite size $J$ of the approximate solutions and their perturbed version within the given classes of perturbations. The related concepts are relevant for then quantify the maximum errors among the real and approximated solutions and parallel issues concerning their counterparts under perturbations of the studied types. More specifically:

**Definition 3.5**. A sampling sequence $\hat{t}_J = \{t_i : i = 0,1,..., J\}$ of strictly ordered sampling points with $h_i = t_{i+1} - t_i \leq h$; $i = 0,1,..., J-1$ is said to be in the class $C_{Jh} = \{t_i \in \hat{t}_J : t_{i+1} - t_i \leq h \; ; \; i = 0,1,..., J-1\}$. □

Note from Definition 3.5 that $h \leq h' \Rightarrow C_{Jh} \subseteq C_{Jh'}$ and that $\hat{t}_J \equiv \{t_i : i = 0,1,..., J-1\} \subset C_{Jh} \Rightarrow t_J - t_0 \leq Jh$.

**Definition 3.6**. A sequence $\hat{x}_J = \{x(t_i) : i = 0,1,..., J-1\}$ of $J$ samples of the solution of an approximate differential equation (2.9) is *a $\delta$-pseudo $J$-orbit of sampling sequence* $\hat{t}_J$ for some $\delta \in R_+$, and denoted by $O(\hat{x}_J, \Gamma, \delta)$ if $\max_{t \in [t_0, t_J]} \|e(t)\| \leq \delta$. □

If the integer $J$ and the real $t_J$ are infinity, the corresponding trajectory solutions are referred to as complete pseudo-orbits and orbits. The solution of the true differential equation (2.4) is a $J$-orbit of sampling sequence $\hat{t}_J$. The continuous approximate (respectively, true) solution for $[t_0, t_J]$ is the $\delta$-pseudo $J$-orbit (respectively, $J$-orbit) of sampling sequence $\hat{t}_J$. The perturbed solutions under the forcing perturbations of Theorem 3.1 and Theorem 3.2 are denoted in a similar way leading to the corresponding perturbed pseudo-orbits.

**Definition 3.7**. The set of all the $\delta$-pseudo-orbits $O(\hat{x}_J, \Gamma, \delta)$ with $\max_{t \in [t_0, t_J]} \|e(t)\| \leq \delta$ for some $\delta \in R_+$ obtained for any sampling sequence $\hat{t}_J \in C_{Jh}$ is said to be *the class $CO(C_{Jh}, \delta)$ of $\delta$-pseudo $J$-orbits of sampling sequence* $\hat{t}_J$. □

**Definition 3.8**. The set $\hat{Y}_J$ of true solution sequences $\hat{y}_J = \{y(t_i) : t_i \in \hat{t}_J, i = 0,1,..., J-1\}$ of sampling sequence $\hat{t}_J$ possesses the shadowing property on the corresponding set of approximate solutions if, for each given $\delta \in R_+$, there is some $y_0 = y_0(\delta)$ for which a $O(\hat{x}_J, \Gamma, \delta)$ exists. It is said that $y_0 = y_0(\delta)$ shadows $O(\hat{x}_J, \Gamma, \delta)$. □

**Proposition 3.9**. If $\hat{Y}_J$ sampling sequence $\hat{t}_J$ possesses the shadowing property then $CO(C_{Jh}, \delta)$ is nonempty for any $\delta \in \mathbf{R}_+$. □

Note that $CO(C_{Jh}, \rho) = \bigcup_{\Gamma \in C_{Jh}} O(\hat{x}_J, \Gamma, \rho)$ and note also that $CO(C_{Jh}, \rho) \subseteq CO(C_{Jh'}, \rho)$ for any $h' \geq h$. The subsequent result relies on Theorem 3.1 and Definition 3.7 for a class of pseudo-orbits $CO(C_{Jh}, \rho)$ defined by a sampling sequence class $C_{Jh}$. In fact, the characterization becomes global for all approximated solutions on a finite interval $[t_0, t_J]$ for sampling intervals $h_i = t_{i+1} - t_i \leq h$; $i = 0, 1, \ldots, J-1$ and initial conditions subject to a maximum allowable deviation with respect to the initial condition of the true solution provided that the approximate solution exists in a global (rather than local) definition domain.

The so-called shadowing properties, [1-5], of the true solutions with respect to the approximated ones has relies on the physical meaning that for sets of appropriate initial conditions, the true solution is arbitrarily close to its approximate version on a certain interval $[t_0, t_J]$ where both solutions exist and are unique. Based on Theorems 2.6, 3.1 and 3.2, the shadowing properties of the true for the approximated solutions, those ones being the nominal one or the perturbed ones under the class of perturbations of Theorems 3.1 and 3.2, are now discussed. It is seen that the shadowing properties at sampling points under Theorems 2.6, 3.1 and 3.2 guarantee the corresponding properties in $[t_0, t_J]$.

The shadowing properties of true solutions of pseudo-orbits for constrained sampling sequences according to the constraints of Theorem 2.6 are addressed in the subsequent result:

**Proposition 3.10**. Consider the true and approximated solutions associated with the differential equations (2.4) and (2.9) satisfying the hypotheses and conditions of Theorem 2.6. Then, such a set of solutions lies in the class $CO(C_{Jh}, \varepsilon)$ of $\varepsilon$-pseudo $J$-orbits of sampling sequence $\hat{t}_J = \hat{t}_J(\rho)$ for $\rho \leq \varepsilon$, subject to one of the constraints (2.12), (2.19)-(2.21) [Lemma 2.4, Theorem 2.6], belonging to a sampling sequence class $C_{Jh}$ for any $\rho, \varepsilon \in \mathbf{R}_+$ with arbitrary $\rho \leq \varepsilon$ and any given $\varepsilon$. Also, there is a $y_0 = y_0(\varepsilon)$ which shadows each $O(\hat{x}_J, \Gamma, \varepsilon) \in CO(C_{Jh}, \varepsilon)$ for each given $\varepsilon \in \mathbf{R}_+$ and $\rho \leq \varepsilon$.

*Proof* One gets from Theorem 2.6 that

$$\max_{t \in [t_0, t_J]} \|e(t)\| \leq \varepsilon = \|e(t_0)\| + \rho$$

for an initial condition $y(t_0)$ of the true differential equation fulfilling $\big| \|y(t_0)\| - \|x(t_0)\| \big| \leq \|e(t_0)\|$ and any given real constants $\varepsilon \geq \rho > 0$. This defines families of initial conditions $y_0 = y_0(\varepsilon)$ of the true differential equation which shadow each $O(\hat{x}_J, \Gamma, \varepsilon) \in CO(C_{Jh}, \varepsilon)$ for each given $\varepsilon \in \mathbf{R}_+$ and $\rho \leq \varepsilon$. For any given $\varepsilon \in \mathbf{R}_+$, it suffices to take $0 < \rho \leq \varepsilon$ to zero in (2.19)-(2.21) of Lemma 2.4 and (2.35)

- (2.36) in Theorem 2.6 to fix an admissible sampling sequence $\hat{t}_J = \hat{t}_J(\varepsilon)$ and then to get the result.
□

The perturbed approximated differential equations referred to in Theorem 3.2 and Theorem 3.1, which is a particular case of Theorem 3.2 for $g(x(t),t)$ being zero for $t \notin \hat{t}_J$, that is for non- sampling points, are analyzed in the subsequent result which generalizes Proposition 3.10:

**Proposition 3.11**. Consider the true and approximated solutions (3.21) and (3.22) associated with the differential equations satisfying the hypotheses and conditions of Theorem 3.2. Then, such a set of solutions lies in the class $CO(C_{Jh}, \varepsilon)$ of $\varepsilon$-pseudo $J$-orbits of sampling sequence $\hat{t}_J = \hat{t}_J(\rho)$ for, subject to one of the constraints (2.12), (2.19) -(2.21), (2.35)-(2.36) [Lemma 2.4, Theorem 2.6] and to either (3.19a) or (3.19b) [Theorem 3.2] with $\rho \leq \varepsilon - \sum_{i=0}^{J-1} \overline{g}_i$ and any arbitrary $\varepsilon \in \mathbf{R}_+$, belonging to a sampling sequence class $C_{Jh}$. Also, there is a $y_0 = y_0(\varepsilon)$ which shadows each $O(\hat{x}_J, \Gamma, \varepsilon) \in CO(C_{Jh}, \varepsilon)$ for each given $\rho, \varepsilon \in \mathbf{R}_+$ with the perturbation fulfilling $\sum_{i=0}^{J-1} \overline{g}_i < \varepsilon$.

*Proof*: One gets from (3.18) in Theorem 3.2 together with either (3.14) or (3.15) that

$$\max_{t \in [t_0, t_J]} \|e(t)\| \leq \varepsilon = \|e(t_0)\| + \rho + \sum_{i=0}^{J-1} \overline{g}_i$$

with any arbitrary real constant $0 < \rho \leq \varepsilon - \sum_{i=0}^{J-1} \overline{g}_i$, provided that $\sum_{i=0}^{J-1} \overline{g}_i < \varepsilon$, and any given real constant $\varepsilon$ for a (shadowing) initial condition $y(t_0)$ of the true differential equation fulfilling

$$\left| \|y(t_0)\| - \|x(t_0)\| \right| \leq \|e(t_0)\| \leq \min\left( \varepsilon - \rho - \sum_{i=0}^{J-1} \overline{g}_i, \frac{1 - \sum_{i=0}^{J-1} h_i \lambda_i}{\sum_{i=0}^{J-1} h_i \lambda_i} \overline{g} \right) \quad (3.32)$$

Note that a sufficient condition guaranteeing (3.31) is

$$\|e(t_0)\| \leq \min\left( \varepsilon - \varepsilon_0 - J\overline{g}, \frac{1 - \sum_{i=0}^{J-1} h_i \lambda_i}{\sum_{i=0}^{J-1} h_i \lambda_i} \overline{g} \right) \quad (3.33)$$

since $\overline{g} \geq \frac{\sum_{i=0}^{J-1} h_i \lambda_i}{1 - \sum_{i=0}^{J-1} h_i \lambda_i} \|e(t_0)\|$. Thus, it suffices to take $0 < \rho \leq \varepsilon$ to zero in either (2.35) - (2.36) in Theorem 2.6 to fix an admissible sampling sequence $\hat{t}_J = \hat{t}_J(\varepsilon)$ and then to get the result. □

**Remark 3.12**. A particular case of Proposition 3.11 for the perturbations (3.4) which are defined only at sampling instants and discussed in Theorem 3.1 is got by the particular constraint below got from (3.32) and (3.33):

$$\|e(t_0)\| \leq \varepsilon - \rho - J\overline{g} \leq \varepsilon - \rho - \sum_{i=0}^{J-1} \overline{g}_i \quad (3.34)$$

□

## ACKNOWLEDGEMENTS

The author is very grateful to the Spanish Government for its support of this research trough Grant DPI2012-30651, and to the Basque Government for its support of this research trough Grants IT378-10 and SAIOTEK S-PE12UN015. He is also grateful to the University of Basque Country for its financial support through Grant UFI 2011/07.